\newenvironment{E}{\begin{equation}}{\end{equation}}
\def\proof{\noindent{\bf Proof: }}
\def\qed{ \hskip 20pt{\vrule height7pt width6pt depth0pt}\hfil}
\def\forb{{\hbox{forb}}}
\def\0{{\bf 0}}
\def\1{{\bf 1}}

\documentclass[12pt]{article}
\renewcommand{\l}{\ell}

\newtheorem{thm}{Theorem}[section]
\newtheorem{lemma}[thm]{Lemma}
\newtheorem{prop}[thm]{Proposition}

\usepackage{amssymb}
\usepackage{amsmath}
\title{Design Theory and Some Non-simple Forbidden Configurations}
\author{R.P.  Anstee\thanks{Research supported in part by
NSERC} and Farzin Barekat\thanks{Research supported by NSERC of first author}
  \\Mathematics Department\\The University of British Columbia\\Vancouver,
B.C. Canada V6T 1Z2\\ {\small\texttt{anstee@math.ubc.ca}}
\\\mbox{\ }}

\begin{document}
\maketitle
\begin{abstract}
Let $\1_k\0_{\l}$ denote the $(k+\l)\times 1$ column of $k$ 1's above $\l$ 0's. Let $q\cdot (\1_k\0_{\l})$ denote the $(k+\l)\times q$ matrix with $q$ copies of the column $\1_k\0_{\l}$. A 2-design $S_{\lambda}(2,3,v)$ can be defined as a $v\times \frac{\lambda}{3}\binom{v}{2}$ (0,1)-matrix with all column sums equal 3 and with no submatrix $(\lambda+1)\cdot(\1_2\0_0)$. Consider an $m\times n$  matrix $A$ with all column sums in $\{3,4,\ldots ,m-1\}$. Assume $m$ is sufficiently large (with respect to $\lambda$) and assume that $A$ has no submatrix which is a row permutation of $(\lambda+1)\cdot (\1_2\0_1)$. Then we show the  number of columns in $A$ is at most $\frac{\lambda}{3}\binom{m}{3}$ with equality for $A$ being the columns of column sum 3 corresponding to the  triples of a 2-design $S_{\lambda}(2,3,m)$. A similar results holds for $(\lambda+1)\cdot(\1_2\0_2)$

Define a matrix to be {\it simple} if it is a (0,1)-matrix with no repeated columns. Given two matrices $A$, $F$, we define $A$ to have $F$ as a {\it configuration} if and only if some submatrix of $A$ is a row and column permutation of $F$.  Given $m$, let $\forb(m,q\cdot (\1_k\0_{\l}))$ denote the maximum number of possible columns in a simple $m$-rowed matrix which has no 
configuration $q\cdot (\1_k\0_{\l})$. For $m$ sufficiently large with respect to $q$, we compute exact values for $\forb(m,q\cdot (\1_1\0_1))$,
$\forb(m,q\cdot (\1_2\0_1))$, $\forb(m,q\cdot (\1_2\0_2))$. In the latter two cases, we use a construction of Dehon (1983) of  {\it simple} triple systems $S_{\lambda}(2,3,v)$ for $\lambda>1$. Moreover for $\l=1,2$, simple $m\times\forb(m,q\cdot(\1_2\0_{\l})$ matrices with no configuration $q\cdot(\1_2\0_{\l})$ must arise from simple 2-designs $S_{\lambda}(2,3,m)$ of appropriate $\lambda$. 

The proofs derive a basic upper bound by a pigeonhole argument and then use careful counting and Tur\'an's bound, for large $m$, to reduce the bound. For small $m$, the larger pigeonhole bounds are sometimes  the exact bound. There are intermediate values of $m$ for which we do not know the exact bound. 
\end{abstract}
\section{Introduction}

Some combinatorial objects can be defined by forbidden substructures. It is also true that most combinatorial objects can be encoded by a (0,1)-matrix. In this paper we consider  submatrices of (0,1)-matrices as the substructures of interest.

Let $\1_k\0_{\l}$ denote the $(k+\l)\times 1$ column consisting of $k$ 1's atop $\l$ 0's. For any positive integer $q$, let $q\cdot (\1_k\0_{\l})$ denote the 
$q\times (k+\l)$ matrix of $q$ copies of $\1_k\0_{\l}$. A 2-design $S_{\lambda}(2,3,v)$ consists of $\frac{\lambda}{3}\binom{v}{2}$ triples from
$\{1,2,\ldots, v\}$ such that for each pair $i,j\in\{1,2,\ldots ,v\}$, there are exactly $\lambda$ triples containing $i,j$. If we encode the triple system as a $v$-rowed (0,1)-matrix $A$ such that the columns are the incidence vectors of the triples, then $A$ has no submatrix $(\lambda+1)\cdot (\1_{2}\0_0)$.
In fact, if $A$ is a $v\times n$ (0,1)-matrix with column sums 3 and $A$ has no submatrix $(\lambda+1)\cdot (\1_2\0_0)$ then
$n\le \frac{\lambda}{3}\binom{m}{2}$ with equality if and only if the columns of $A$ correspond to the triples of a 2-design $S_{\lambda}(2,3,v)$. This can be shown by a pigeonhole counting argument. 

The problem of forbidding a submatrix is usually extended to forbidding any row and column permutation of the submatrix. Let $A$, $F$ be (0,1)-matrices. We say that $A$ has $F$ as a {\it configuration} if there is a submatrix of $A$ which is a row and column permutation of $F$. We extend the forbidden submatrix $(\lambda+1)\cdot(\1_2\0_0)$ and obtain the following two design theory results.

\begin{thm}\label{designq110} Let $\lambda$ and $v$ be given integers. There exists an $M$ so that for $v>M$, if $A$ is an $v\times n$ (0,1)-matrix with column sums in $\{3,4,\ldots ,v-1\}$and  $A$ has no configuration $(\lambda+1)\cdot(\1_2\0_1)$ then
\begin{E}n\le \frac{\lambda}{3}\binom{v}{2}\end{E}
and we have equality if and only if the columns of $A$ correspond to the triples of a
2-design $S_{\lambda}(2,3,v)$.\qed\end{thm}

When we extend the forbidden configuration to $(\lambda+1)\cdot(\1_2\0_2)$ the case of equality becomes more difficult.

\begin{thm}\label{designq1100} Let $\lambda$ and $v$ be given integers. There exists an $M$ so that for $v>M$, if $A$ is an $v\times n$ (0,1)-matrix with column sums in $\{3,4,\ldots ,v-3\}$ and  $A$ has no configuration $(\lambda+1)\cdot(\1_2\0_2)$ then
\begin{E}n\le \frac{\lambda}{3}\binom{v}{2}\end{E}
and we have equality if and only if there are positive integers $a,b$ satisfying $a+b=\lambda$ and there are $\frac{a}{3}\binom{v}{2}$ columns of $A$ of column sum 3 corresponding to the triples of a
2-design $S_{a}(2,3,v)$  and there are $\frac{b}{3}\binom{v}{2}$ columns of $A$ of column sum $v-3$ of $v-3$-sets whose complements (in $\{1,2,\ldots ,v\}$) corresponding to the triples of a
2-design $S_{b}(2,3,v)$.\qed\end{thm}

Our first motivation for studying these problems came from 
extremal set theory. 
An $m\times n$ (0,1)-matrix $A$ can be thought of a multiset of $n$
subsets of $\{1,2,\ldots ,m\}$. Let $[m]=\{1,2,\ldots ,m\}$. For an $m\times 1$ (0,1)-column $\alpha$, we define 
\begin{E}S(\alpha)=\{i\in [m]\,:\, \alpha\hbox{ has }1\hbox{ in row }i\}.\label{salpha}\end{E}
From this we define the natural multiset system $\cal A$ associated with the matrix $A$:
\begin{E}{\cal A}=\{S(\alpha_i)\,:\, \alpha_i\hbox{ is  column }i\hbox{ of }A \}.\label{calA}\end{E}
Similarly, if we are given a multiset system ${\cal A}$, we can form a matrix $A$, as long as we don't care about column order. We define a {\it simple} matrix $A$ as a (0,1)-matrix with no repeated columns. In this case ${\cal A}$ yields as {\it set} system and it is in this setting that extremal set theory problems can be stated.

We define  $\forb(m,F)$ as the smallest value (depending on $m$ and $F$) so
that if $A$ is a simple $m\times n$ matrix and  $A$ has no
configuration $F$ then $n\le \forb(m,F)$. Alternatively $\forb(m,F)$
is the smallest value so that if $A$ is an $m\times
(\forb(m,F)+1)$ simple  matrix then $A$ must have a configuration $F$. 
A sampling of exact results for $\forb(m,F)$ are in \cite{AFS}, \cite{AKa}.  

Let $K_k$ denote the
$k\times 2^k$ simple matrix of all possible (0,1)-columns on $k$ rows and let
$K_k^s$ denote the $k\times \binom{k}{s}$ simple matrix of all
possible columns of column sum $s$. Many
results have been obtained about $\forb(m,F)$. Exact results have been rare for non-simple configurations  $F$. We consider $F=q\cdot(\1_k\0_{\l}$) for $(k,\l)=(1,1),(2,1),(2,2)$. In \cite{AFS} we showed that 
$$\left\lfloor\frac{q+1}{2}m\right\rfloor+2\le \forb(m, q\cdot(\1_1\0_1))\le
\left\lfloor\frac{q+1}{2}m+\frac{(q-3)m}{2(m-2)}\right\rfloor+2$$
    where the upper bound obtained by a pigeonhole argument is achieved for $m=q-1$ by taking $A=[K_m^0K_m^1K_m^2K_m^{m-1}K_m^m]$.  For $m$ with $m\ge \max\{3q+2, 8q-19\}$, we are able to show that the lower bound is correct and  slice $\frac{(q-3)m}{2(m-2)}\approx \frac{q-3}{2}$ off the pigeonhole bound. It is likely that our bound is valid for smaller $m>q-1$. The case $q=4$, is Lemma 3.1 in \cite{AKa} and took a page to establish.

\begin{thm}\label{q10} Let $q\ge 3$ be given. Then for $m\ge \max\{3q+2, 8q-19\}$, \begin{E}\forb(m, q\cdot (\1_1\0_1)=\Bigl[\left.
\overbrace{\begin{matrix}1& 1& \cdots &1\\ 0 &0& \cdots &0\end{matrix}}^{q}\right.\Bigr])=\lfloor\frac{q+1}{2}m\rfloor +2.
\hskip 1in\qed\label{q10bd}\end{E}
\end{thm}

For $m$ even or $q-3$ even, let $G$ be a (simple) graph on $m$ vertices for which all the degrees are $q-3$ and for $m,q-3$ odd let $G$ be a graph for which $m-1$ vertices have degree $q-3$ and one vertex has degree $q-4$. Such graphs are easy to construct.
Let $H$ be the vertex-edge incidence matrix associated with $G$, namely for each edge $e=(i,j)$ of $G$, we add a column to $H$ with 1's in rows $i,j$ and 0's in other rows. Thus $H$ is a simple $m$-rowed matrix  with 
$\lfloor\frac{(q-3)m}{2}\rfloor$ columns each of column sum 2.
The simple matrix $A=[K_m^0\,K_m^1\,H\,K_m^{m-1}\,K_m^m]$ has $\lfloor\frac{(q+1)m}{2}\rfloor+2$ columns  and no configuration
$q\cdot (\1_1\0_1)$ which establishes $\forb(m,q\cdot (\1_1\0_1))\ge \lfloor\frac{(q+1)m}{2}\rfloor+2$. We establish the upper bound in Section~\ref{sectionq10}. 

We are able to solve two more cases but need certain designs to achieve exact bounds.  A 2-design $S_{\lambda}(2,3,v)$ (or {\it triple system}) is defined to be {\it simple} if no triple is repeated. The associated $v\times\frac{\lambda}{3}\binom{v}{2}$ matrix is a simple matrix. We need the following result.

\begin{thm}\label{design}Dehon\cite{D} Let $v,\lambda$ be given. Then a simple     2-design $S_{\lambda}(2,3,v)$ exists if and only if $v(v-1)\equiv 0(\hbox{mod }6)$, $v-1\equiv 0(\hbox{mod }2)$ and $v\ge \lambda+2$. \qed\end{thm}

These designs are used in the constructions for the following two theorems in the following way. We form a simple $v\times \frac{\lambda}{3}\binom{v}{2}$  matrix  $T_{v,\lambda}$  whose columns correspond to the blocks of $S_{\lambda}(2,3,v)$  so that if $B$ is a block then the corresponding column has a $1$ in row $i$ if and only if $i\in B$. Note that $T_{v,\lambda}$  has no submatrix 
$(\lambda+1)\cdot (\1_2\0_0)$. Pigeonhole arguments will show  that
$\forb(m,q\cdot (\1_2\0_0))\le \binom{m}{0}+\binom{m}{1}+\frac{q+1}{3}\binom{m}{2}$ with equality, by Dehon's Theorem~\ref{design}, for $m\ge q$ and $m\equiv 1,3(\hbox{mod }6)$.  The matrix achieving equality would be $[K_m^0K_m^1K_m^2T_{m,q-2}]$. Let $B^c$ denote the (0,1)-complement of a matrix $B$. Note that  the $v\times \frac{a+b}{3}\binom{v}{2}$ simple matrix 
$[T_{v,a}\,T_{v,b}^c]$ has no submatrix 
$(a+b+1)\cdot (\1_2\0_0)$.

\begin{thm}\label{q110} Let $q>2$ be given. There exists a constant $M=M(q)$ so that for $m>M$, 
\begin{E}\forb(m,q\cdot (\1_2\0_1)=
\left[\begin{matrix}\\ \\ \\\end{matrix}\right.
\overbrace{\begin{matrix}1& 1& \cdots &1\\ 1&1&\cdots &  1\\ 0 &0& \cdots &0\end{matrix}}^{q}
\left.\begin{matrix}\\ \\ \\\end{matrix}\right]
)\le m+2+\frac{q+1}{3}\binom{m}{2}\label{simpleq110bd}\end{E}
with equality for $m\equiv 1,3(\hbox{mod }6)$. If
$A$ is an $m\times\forb(m,q\cdot (\1_2\0_2))$ simple matrix with $m>M$ and $m\equiv 1,3(\hbox{mod }6)$, then $A$ consists of all possible columns of sum 0, 1, 2, $m$ and the columns of column sum 3 correspond to a simple triple system $T_{m,q-2}$ and $A$ has no further columns. \qed\end{thm}

\begin{thm}\label{q1100} Let $q>2$ be given. There exists a constant $M=M(q)$ so that for $m>M$, 
\begin{E}\hbox{forb}(m, q\cdot (\1_2\0_2)=
\left[\begin{matrix}\\ \\ \\ \\ \end{matrix}\right.
\overbrace{\begin{matrix}1& 1& \cdots &1\\ 1&1&\cdots &  1\\ 0 &0& \cdots &0\\ 0 &0& \cdots &0\\ \end{matrix}}^{q}
\left.\begin{matrix}\\ \\ \\ \\ \end{matrix}\right])\le 2+2m+\frac{q+3}{3}\binom{m}{2},\label{simpleq1100bd}\end{E}
with equality for $m\equiv 1,3(\hbox{mod }6)$. 
If $A$ is an $m\times\forb(m,q\cdot (\1_2\0_2))$ simple matrix with $m>M$ and $m\equiv 1,3(\hbox{mod }6)$,  then there exist positive integers $a,b$ with $a+b=q-3$ so that $A$ consists of all possible columns of sum 0, 1, 2, $m-2$, $m-1$, $m$ and the columns of column sum 3 correspond to a simple triple system $T_{m,a}$ and the columns of column sum $m-3$ correspond to the complement of a simple triple system $T_{m,b}$ and $A$ has no further columns.
\qed\end{thm}

Thus the constructions for equality in Theorem~\ref{q110} are $A=[K_m^0K_m^1K_m^2T_{m,q-2}K_m^m]$ and the constructions  for equality in Theorem~\ref{q1100} are  found by selecting $a,b$ positive integers with $a+b=q-3$ and using
 $A=[K_m^0K_m^1K_m^2T_{m,a}T_{m,b}^cK_m^{m-2}K_m^{m-1}K_m^m]$. For $m=q+1$, the construction
$A=[K_m^0K_m^1K_m^2K_m^3K_m^m]$ avoids $q\cdot(\1_2\0_1)$ and exceeds the bound (\ref{q110bd}) and for $m=q+1$, the construction
$A=[K_m^0K_m^1K_m^2K_m^3K_m^{m-2}K_m^{m-1}K_m^m]$ and avoids $q\cdot(\1_2\0_2)$ and exceeds the bound (\ref{q1100bd}) so our theorems need some condition on $m$.

To prove Theorem~\ref{designq110} and Theorem~\ref{q110}, we prove the following:
 
\begin{prop}\label{nonsimpleq110} Let $m$, $q>2$ be given. Let $A$ be an $m\times n$ (0,1)-matrix so that no column of sum 0,1,2, or $m$ is repeated. Assume $A$ has no configuration $q\cdot (\1_2\0_1)$. Then there exists a constant $M$ so that for $m>M$, 
\begin{E}n\le
  m+2+\frac{q+1}{3}\binom{m}{2}\label{q110bd}\end{E}
with equality for $m\equiv 1,3(\hbox{mod }6)$.
If $A$ is an $m\times\forb(m,q\cdot (\1_2\0_1))$ matrix with $m>M$ and $m\equiv 1,3(\hbox{mod }6)$, then $A$ consists of all possible columns of sum 0, 1, 2, $m$ once each and the columns of column sum 3 correspond to the triples of a 2-design $S_{q-2}(2,3,m)$ and $A$ has no further columns. \qed\end{prop}

We see that Theorem~\ref{designq110} follows by taking a matrix $A$ of column sums in $\{3,4,\ldots ,m-1\}$ and with no configuration $(\lambda+1)\cdot(\1_2\0_1)$ and adding the $\binom{m}{2}+m+2$ columns of column sum 0,1,2 and $m$ to obtain a matrix $A'$. Now $A'$ has no configuration $(\lambda+2)\cdot (\1_2\0_1)$ and satisfies the hypotheses of Proposition~\ref{nonsimpleq110} with $q=\lambda+2$. Applying Proposition~\ref{nonsimpleq110} yields Theorem~\ref{designq110}. The bound of Theorem~\ref{q110}
follows  directly from Proposition~\ref{nonsimpleq110}. 
To prove Theorem~\ref{designq1100} and Theorem~\ref{q1100} we prove the following:
 
\begin{prop}\label{nonsimpleq1100} Let $m$, $q>2$ be given. Let $A$ be an $m\times n$ (0,1)-matrix so that no column of sum 0,1,2, $m-2$, $m-1$ or $m$ is repeated. Assume $A$ has no configuration $q\cdot (\1_2\0_2)$. Then there exists a constant $M$ so that for $m>M$, 
\begin{E}n\le
  2m+2+\frac{q+3}{3}\binom{m}{2}\label{q1100bd}\end{E}
with equality for $m\equiv 1,3(\hbox{mod }6)$.
If $A$ is an $m\times\forb(m,q\cdot (\1_2\0_2))$ matrix with $m>M$ and $m\equiv 1,3(\hbox{mod }6)$, then $A$ consists of all possible columns of sum 0, 1, 2, $m-2$, $m-1$ and $m$ once each and there are two positive integers $a,b$ satisfying $a+b=q-3$ with the columns of column sum 3 correspond to the triples of a 2-design $S_{a}(2,3,m)$ and the columns of column sum $m-3$ correspond to the complements in $[m]$ of the blocks of a 2-design $S_{b}(2,3,m)$ and $A$ has no further columns. \qed\end{prop}

We see that Theorem~\ref{designq1100} follows by taking a matrix $A$ of column sums in $\{3,4,\ldots ,m-3\}$ and with no configuration $(\lambda+1)\cdot(\1_2\0_2)$ and adding the $2\binom{m}{2}+2m+2$ columns of column sum 0,1,2, $m-2$, $m-1$ and $m$ to obtain a matrix $A'$. Now $A'$ has no configuration $(\lambda+3)\cdot (\1_2\0_2)$ and satisfies the hypotheses of Proposition~\ref{nonsimpleq1100} with $q=\lambda+3$. Applying Proposition~\ref{nonsimpleq1100} yields Theorem~\ref{designq1100}. The bound of Theorem~\ref{q1100}
follows  directly from Proposition~\ref{nonsimpleq1100}.

We could give a simpler direct proof of Theorem~\ref{designq110} by using the proof of Proposition~\ref{nonsimpleq1100} and  ignoring certain column sums. We were originally motivated by  the forbidden configuration bounds of  Theorems~\ref{q110} and Theorem~\ref{q1100}.
 
The proofs of Proposition~\ref{nonsimpleq110} and Proposition~\ref{nonsimpleq1100} use Tur\'an's bound for the maximum number of edges in a graph with no complete graph of a certain size. We do not explicitly give values for $M$ since the values as given by the proofs are unlikely to be of value but our proof shows we may take $M$ to be $O(q^3)$. Proposition~\ref{nonsimpleq110} for $q\cdot(\1_2\0_1)$ is proven in Section~\ref{sectionq110} and Proposition~\ref{nonsimpleq1100} for $q\cdot(\1_2\0_2)$ is proven in Section~\ref{sectionq1100}. The proofs are organized to highlight analogies with the proof of Theorem~\ref{q10} but the details are different.
We were surprised that  exact bounds were obtained.
We do not see how to extend our exact proofs to $F=t\cdot (\1_k\0_k)$ with $k\ge 3$ and moreover do not have the analogue of Dehon's lovely Theorem~\ref{design} to provide a construction of simple $k$-designs.

\section{Exact Bound for $q\cdot(\1_1\0_1)$}\label{sectionq10}

This section gives the proof of Theorem~\ref{q10}. We have broken it into  lemmas.
Assume $A$ is a simple $m$-rowed matrix with no configuration $q\cdot(\1_1\0_1)$.
Let $a_i$ denote the number of columns with either exactly $i$ 1's or $i$ 0's for $i=0,1,2$ and let $a_3$ be the number of remaining columns. Without loss of generality, we may assume $a_0=2$ since the column of all 0's and the column of all 1's cannot contribute to $q\cdot (\1_1\0_1)$. Thus $2+a_1+a_2+a_3$ is the number of columns of $A$.

In \cite{AFS}, we establish that
$$2+a_1+a_2+a_3\le 
\left\lfloor\frac{(q+1)m}{2}+\frac{(q-3)m}{2(m-2)} \right\rfloor+2$$
and  as noted in the Introduction, we can achieve equality for some small $m$. We wish to show that these small values of $m$ are exceptional. We assume
\begin{E}a_1+a_2+a_3>\left\lfloor\frac{(q+1)m}{2} \right\rfloor\label{a1a2a3}\end{E}
and seek a contradiction.

\begin{lemma} Let $A$ be an $m\times n$ simple matrix with no $q\cdot (\1_1\0_1)$. Assume $m\ge 6$. Then
\begin{E}(m-1)a_1+2(m-2)a_2+ 3(m-3)a_3\le (2q-2)\binom{m}{2}
=(q-1)m(m-1).\label{risefallineq}\end{E}
Assume  $n> \frac{q+1}{2}m+2$. Then
\begin{E}2m-\frac{m(q-3)}{m-3}<a_1\le 2m,\label{a1ineq}\end{E}
\begin{E}a_3<\frac{m(q-3)}{m-5}.\label{a3ineq}\end{E}
\label{pigeonholeq10}\end{lemma}

\proof A column of $k$ 1's contains $\binom{k}1\binom{m-k}1$ configurations $\1_1\0_1$. Note that $\binom{k}{1}\binom{m-k}{1}\ge \binom{3}{1}\binom{m-3}{1}$ for $3\le k\le m-3$. By the pigeonhole argument, there are at most $(2q-2)\binom{m}{2}$ configurations $\1_1\0_1$ in $A$ else there will be $2q-1$ in one of the $\binom{m}{2}$ pairs of rows and hence at least $q$ with the 1 of the $\1_1\0_1$ in the same row yielding the configuration $q\cdot (\1_1\0_1)$. This yields (\ref{risefallineq}).
Given $m\ge 6$, we have $m-1<2(m-2)<3(m-3)$. Substituting in (\ref{risefallineq}),

$$a_1(m-1)+2(m-2)(a_2+a_3)\le m(m-1)(q-1).$$
Using $a_2+a_3>\frac{q+1}{2}m-a_1$ from (\ref{a1a2a3})  we have

$$ a_1(m-1)+2(m-2)( \frac{q+1}{2}m-a_1)< m(m-1)(q-1)$$ 
and so
$$2m^2-mq-3m<(m-3)a_1$$
from which we deduce the lower bound of (\ref{a1ineq}). The upper bound of
(\ref{a1ineq}) follows from counting all possible columns.

To show $a_3$ is small, use (\ref{risefallineq}) to obtain
$$a_1(m-1)+2(m-2)(\frac{q+1}{2}m-a_1-a_3)+3(m-3)a_3< m(m-1)(q-1).$$
Rearranging yields
$$(m-5)a_3<m(q-2m+3)+(m-3)a_1.$$
Substituting $a_1\le 2m$, we obtain (\ref{a3ineq}). \qed

\vskip 5pt
Form two graphs $G_0$, $G_1$ from the columns of $A$ where the vertex set for both graphs corresponds to the  rows of $A$. We form a graph $G_0$ from the columns of $A$ of column sum $m-2$ so that if there is a column of $A$ with $m-2$ 1's and two $0$'s on rows $i,j$ we add an edge $(i,j)$ to $G_0$. Similarly we form a graph $G_1$ from the columns of $A$ of column sum 2 so that if there is a column of $A$ with $m-2$ 0's and two $1$'s on rows $i,j$, then $G_1$ has the edge $(i,j)$. Define $d_0(i)$ and $d_1(i)$ to be the degrees of $i$ in $G_0$ and $G_1$ respectively. Hence
\begin{E}a_2=\frac{1}{2}\sum_{i=1}^m(d_0(i)+d_1(i)).\label{a2eq10}\end{E}
Using (\ref{a1a2a3}),  we obtain
$$a_1+\frac{1}{2}\sum_{i=1}^m(d_0(i)+d_1(i))+a_3>\frac{q+1}{2}m$$
Multiplying by 2 and substituting the upper bounds (\ref{a1ineq}) for $a_1$ and (\ref{a3ineq}) for $a_3$, yields
$$\sum_{i=1}^m(d_0(i)+d_1(i))>(q+1)m-4m-\frac{2m(q-3)}{m-5}$$
\begin{E}=m(q-3)\Bigl(1-\frac{2}{m-5}\Bigr).\label{2a2ineq}\end{E}
Thus the average value of $d_0(i)+d_1(i)$ is close to $q-3$.

The possible columns of column sum 1 or $m-1$ are as follows.
Define $e_i$ to be the $m$-rowed column with a 1 in row $i$ and 0's elsewhere and let $e_i^c$ be the (0,1)-complement of $e_i$. Define
$$E_1=\{i\,:\,1\le i\le m \hbox{ and }e_i\hbox{ is not in $A$}\},$$ 
$$E_0=\{i\,:\,1\le i\le m \hbox{ and }e_i^c\hbox{ is not in $A$}\}.$$

We have $a_1=2m-|E_0|-|E_1|$ and so $|E_1|+|E_0|<\frac{m(q-3)}{m-3}$ by (\ref{a1ineq}). For convenience of counting define
\begin{E}\epsilon(i)=\left\{\begin{array}{ll}0&\hbox{ if }i\notin E_1\cup E_0\\
1&\hbox{ if }i\in E_1\backslash E_0\hbox{ or }i\in E_0\backslash E_1\\ 
2&\hbox{ if }i\in E_1\cap E_0\\\end{array}\right.\quad .\label{e1e0}\end{E}
Thus $\sum_{i=1}^m\epsilon(i)= |E_0|+|E_1| $.

\begin{lemma} Assume $m>3q+2$. Then for all $i=1,2,\ldots ,m$, we have $d_0(i)+d_1(i)\le q-3+\epsilon(i)$. \label{dupperbd}\end{lemma}
\proof
Assume the contrary that  $k$ is an index with $\l=d_0(k)+d_1(k)
\ge q-2+\epsilon(k)$.
Let $N_1$ be the vertices/rows connected to $k$ by no edges in either $G_0$ or $G_1$. Let $N_2$ be the number of vertices connected to $k$ by an edge in $G_0$ or an edge in $G_1$ but not both. Let $N_3$ be the number of vertices connected to $k$ by edges in both $G_0$ and $G_1$. We have
\begin{E}|N_1|+|N_2|+|N_3|=m-1,\qquad |N_2|+2|N_3|= d_0(k)+d_1(k)=\l .\label{n22n3}\end{E}

Consider a row $i\ne k$. There are
at most $2q-2$ configurations $\1_1\0_1$ contained in rows $k,i$ of $A$ and there are $4-\epsilon(i)-\epsilon(k)$ configurations $\1_1\0_1$ contained in rows $k,i$ of $A$ in the columns of column sum $1$ or $m-1$ (corresponding to those columns $e_k,e_k^c,e_i,e_i^c$ which are present in $A$).
If $i\in N_1$ then each edge incident with either $k$ or $i$ in either $G_0$ or $G_1$ corresponds to a column of $A$ of column sum $2$ or $m-2$ that has the configuration $\1_1\0_1$ in rows $k,i$ and hence we have 
$ d_1(k)+d_0(k)+d_1(i)+d_0(i)$ configurations  $\1_1\0_1$ in these columns.
Thus 
$d_1(k)+d_0(k)+d_1(i)+d_0(i)+ (4-\epsilon(i)-\epsilon(k) )\le 2q-2 $ which yields
$$d_1(i)+d_0(i)\le 2q-6-\l+ \epsilon(i)+\epsilon(k). $$

In the case $i\in N_2$ then we note that an edge in say $G_0$ joining $k,i$ contributes 2 to $d_0(i)+d_0(k)$ but the corresponding column does not contain the configuration $\1_1\0_1$ in rows $i,k$. A similar argument holds for an edge $(k,i)$ in $G_1$. By the above analysis we obtain
$$d_1(i)+d_0(i)\le 2q-4-\l+ \epsilon(i)+\epsilon(k). $$

In the case $i\in N_3$ then we note that the two edges in $G_0$ and $G_1$ joining $k,i$ contributes 4 to $d_0(i)+d_1(i)+d_0(k)+d_1(k)$ but correspond to only two columns neither of which  contain the configuration $\1_1\0_1$. By the above analysis we obtain
$$d_1(i)+d_0(i)\le 2q-2-\l+ \epsilon(i)+\epsilon(k). $$
Summarizing, we have for $i\in N_j$ and $j=1,2,3$ that
\begin{E}d_1(i)+d_0(i)\le 2q-6+2(j-1)-\l+ \epsilon(i)+\epsilon(k). \label{dbound}\end{E}

Now we sum our upper bounds on $d_0(i)+d_1(i)$ over all rows $i\in[m]=\{k\}\cup N_1\cup N_2\cup N_3$ and use (\ref{2a2ineq}) to obtain
$$\l+ \sum_{j\in\{1,2,3\}}\sum_{i\in N_j}\bigl(2q-6+2(j-1)-\l+\epsilon(i)+\epsilon(k)\bigr) 
$$
$$\ge \sum_{i=1}^m(d_0(i)+d_1(i)) >m(q-3)\biggl(1-\frac{2}{m-5}\biggr)$$
This simplifies to
$$\l+(2q-6)(|N_1|+|N_2|+|N_3|)+2(|N_2|+2|N_3|)- (m-1)\l+$$
$$+(|E_0|+|E_1| - \epsilon(k))+(m-1)\epsilon(k)
 >m(q-3)(1-\frac{2}{m-5})$$
Using $|N_1|+|N_2|+|N_3|=m-1$, $\l=|N_2|+2|N_3|$ from (\ref{n22n3}), and $|E_0|+|E_1| \le \frac{m(q-3)}{m-3}$  and rearranging yields
$$(2q-6)(m-1)-(m-4)\l+(m-4)\epsilon(k)+2\epsilon(k)+ \frac{m(q-3)}{m-3}
>m(q-3)-\frac{2m(q-3)}{m-5}$$
Using $-\l+\epsilon(k)\le -(q-2)$ and $\epsilon(k)\le 2$ and rearranging we get
\begin{E}\frac{m(q-3)}{m-3}+\frac{2m(q-3)}{m-5}>m-2.\label{messy}\end{E}
We can rewrite (\ref{messy}) as
$0>m^3-(3q+2)m^2+(11q-2)m-30$
which is impossible for $m> 3q+2$. This contradiction establishes the lemma. \qed

\vskip 5pt
Let 
$$Y=\{i\,:\,d_0(i)+d_1(i)=q-3\hbox{ and }\epsilon(i)=0\}$$

\begin{lemma} Assume $m>\max\{3q+2,8q-19\}$.  
Then we may assume $|Y|\ge m/2$.
\end{lemma}
\proof We consider $[m]$ divided into $Y$, $E_0\cup E_1$, and $[m]\backslash (Y\cup E_0\cup E_1)$.
We use Lemma~\ref{dupperbd}.  We have
$$\sum_{i\in E_0\cup E_1}d_0(i)+d_1(i)\le \sum_{i\in E_0\cup E_1}(q-3)+ |E_0|+|E_1|  = |E_0\cup E_1|(q-3)+ |E_0|+|E_1|$$
using $\sum_{i=1}^m\epsilon(i)= |E_0|+|E_1|$. We readily compute
$\sum_{i\in Y}d_0(i)+d_1(i)= |Y|(q-3)$ and
$$ \sum_{i\in [m]\backslash (Y\cup E_0\cup E_1)}d_0(i)+d_1(i)\le 
\sum_{i\in [m]\backslash (Y\cup E_0\cup E_1)} (q-4)\le (m-|Y|-|E_0\cup E_1|)(q-4) $$
Summing we obtain
$$\sum_{i\in[m]}d_0(i)+d_1(i)\le m(q-3)+|E_0|+|E_1|-m+|Y|+|E_0\cup E_1|$$
Now using (\ref{2a2ineq}),  we deduce
$$|E_0|+|E_1|+|E_0\cup E_1|+\frac{2m(q-3)}{m-5}> m-|Y|$$
We use $|E_0\cup E_1|\le|E_0|+|E_1|<\frac{m(q-3)}{m-3}$ by (\ref{a1ineq}) to obtain
$\frac{2m(q-3)}{m-3}+\frac{2m(q-3)}{m-5}> m-|Y|$.  Now for
$m>8q-19$ (so that $m-3>m-5\ge 8(q-3)$), we have 
$\frac{2m(q-3)}{m-3}+\frac{2m(q-3)}{m-5}\le m/2$.
Thus for $m>8q-19$, we may assume $|Y| \ge m/2$.
\qed
\vskip 5pt
Let $A_3$ denote the submatrix of $A$ formed by the columns of sum $3,4,\ldots ,$ or  $m-3$. Then $A_3$ has $a_3$ columns. Let $A_3(Y)$ denote the submatrix of $A_3$ indexed by the rows of $Y$.

\begin{lemma} Assume $m>\max\{3q+2,8q-19\}$. Then $A_3(Y)$ has no configuration $\1_1\0_1$. \label{no10}\end{lemma}

\proof Assume there is a column $\alpha$ in $A_3$ which has both 0's and 1's in the rows indexed by $Y$. By taking the (0,1)-complement of $A$ if necessary, we may assume the number of 1's in those rows is at least $|Y|/2\ge m/4$.  Consider   a row $i\in Y$ where $\alpha$ has a 0. Then there exists a row $j\in Y$ where $\alpha$ has a 1 such that rows $i,j$ are not connected in $G_0$ or $G_1$, since $i$ is connected to at most $q-3$ rows and $|Y|/2\ge m/4>q-3$. 
Given $i,j\in Y$, we have $d_0(i)+d_1(i)=d_0(j)+d_1(j)=q-3$ and $\epsilon(i)=\epsilon(j)=0$.
Given that $i,j$ are not connected in $G_0$ or $G_1$, we have $2(q-3)$ copies of the configuration $\1_1\0_1$ on rows $i,j$ in the columns (of $A$) of exactly two 1's or exactly two 0's. Given $Y\cap (E_1\cup E_2)=\emptyset$, we have $4$ copies of the configuration $\1_1\0_1$ in the columns (of $A$) of one 1 or one 0 and in rows $i,j$.  But $\alpha$ has $\1_1\0_1$ in rows $i,j$ and so we find $q\cdot(\1_1\0_1)$ in $A$, a contradiction. This establishes the lemma. \qed

\vskip 5pt
\noindent{\bf Proof of Theorem~\ref{q10}}: We obtain a contradiction from assuming (\ref{a1a2a3}) and $m>\max\{3q+2,8q-19\}$ and thus establish (\ref{q10bd}). By Lemma~\ref{no10}, each column of $A_3$ has either all 1's or all 0's on the rows of $Y$. Considering column sums, every column in $A_3$ which has all 1's on rows $Y$, has at least three 0's and every column in $A_3$ which has all 0's on rows $Y$, has at least three 1's. 
For $i\in[m]\backslash Y$, let $t_0(i)$ denote the number of 0's in columns of column sum in 
$\{3,4,\ldots ,m-3\}$ which are all 1's on $Y$ and let $t_1(i)$ denote the number of 1's in columns of column sum in 
$\{3,4,\ldots ,m-3\}$ which are all 0's on $Y$. Counting yields 
\begin{E}\sum_{i\in [m]\backslash Y}(t_0(i)+t_1(i))\ge 3a_3.\label{tbound}\end{E}
Let $i\in [m]\backslash Y$ be given. We wish to establish
\begin{E}d_0(i)+d_1(i)\le q-3+\epsilon(i)-t_0(i)-t_1(i)
\label{dupperbdwitht}\end{E}
We  use a similar argument as Lemma~\ref{dupperbd}. Consider a column of $A_3$ which is all 0's on rows of $Y$. Then the column has a configuration $\1_1\0_1$ in rows $i,k$ for any choice of $k\in Y$. A similar remarks holds for columns of $A_3$ which are all 1's on rows of $Y$. Let $X$ denote all the neigbours of $i$ in $G_0$ and in $G_1$. We have $|X|\le d_0(i)+d_1(i)\le q-1$ using Lemma~\ref{dupperbd}. Given $|Y|>m/2>q-1$, we can select a $k\in Y$ with $k\notin X$. Now the columns of sum 1 or $m-1$ in $A$ have
$4-\epsilon(i)-\epsilon(k)= 4-\epsilon(i)$ configurations $\1_1\0_1$ in rows $i,k$ (since $\epsilon(k)=0$). Given that $k\notin X$, the columns of column sum $2$ or $m-2$ have $d_0(i)+d_1(i)+d_0(k)+d_1(k)=d_0(i)+d_1(i)+q-3$ configurations
$\1_1\0_1$ in rows $i,k$. The columns of sum at least 3 and at most $m-3$ have at least $t_0(i)+t_1(i)$ configurations
$\1_1\0_1$ in rows $i,k$ for that choice of $k$.  Rows $i,k$ of $A$ have at most $2(q-1)$ such configurations and so we obtain (\ref{dupperbdwitht}).

Combining twice (\ref{a1a2a3}) and (\ref{a2eq10}) we have
$$2a_1+\sum_{i=1}^m(d_0(i)+d_1(i))+ 2a_3>m(q+1).$$
Using $a_1=2m-|E_0|+|E_1| $, substituting $d_0(i)+d_1(i)=q-3$ for $i\in Y$ and using (\ref{dupperbdwitht}),
$$\sum_{i\in Y}(q-3)+\sum_{i\in [m]\backslash Y}(q-3+\epsilon(i)-t_0(i)-t_1(i))+2a_3>m(q+1)-2(2m-|E_0|+|E_1|).$$
Now using (\ref{tbound}) and $\sum_{i\in [m]\backslash Y}\epsilon(i)\le |E_0|+|E_1|$,
$$|Y|(q-3)+(m-|Y|)(q-3) +|E_0|+|E_1| -a_3 >m(q-3)+2(|E_0|+|E_1| )$$
which yields the contradiction (even for $a_3=0$ and $|E_0|+|E_1| =0$)
$$-a_3>|E_0|+|E_1|.$$
This final contradiction  establishes (\ref{q10bd}).
\qed

One could note that for a matrix $A$ to achieve equality, we would have $a_3=0$ and $|E_0|+|E_1|=0$ and so $a_1=2m$. This suggests that $A$ would have to correspond to the construction given in the Introduction or its (0,1)-complement.

\section{Exact Bound for $q\cdot(\1_2\0_1)$}\label{sectionq110}

We are able to generalize the argument for Theorem~\ref{q10} following a similar series of Lemmas to obtain
Proposition~\ref{nonsimpleq110}. We do not explicitly calculate the smallest possible constant $M$  for our proof (following the argument yields that $M$ is $O(q^3)$), believing that our argument does not give a realistic values for $M$.
Let $A$ be an $m\times n$ (0,1)-matrix with no configuration $q\cdot(\1_2\0_1)$ so that there are no repeated columns of sum $0,1,2,m$. We wish to ignore the $m+2$ possible columns of sum 0, 1, $m$ since they cannot contribute to a configuration $q\cdot (\1_2\0_1)$. So assume $A$ has column sums between 2 and $m-1$, inclusive. Assume $n>\frac{q+1}{3}\binom{m}{2}$. We wish to arrive at a contradiction to prove (\ref{q110bd}).

For $i=2,3$, let $a_i$ denote the number of columns of column sum $i$ in $A$ and let $a_4$ denote the number of columns of column sum at least 4 in $A$.  Note that the  definition of $a_1,a_2,\ldots$ is different in this section from Sections~\ref{sectionq10} and \ref{sectionq1100}. Note that we do not allow repeated columns of sum 2.
We have by assumption that
\begin{E}a_2+a_3+a_4> \frac{q+1}{3}\binom{m}{2}.\label{a2a3a4110}\end{E}
\begin{lemma}\label{pigeon110} Let $m,q$ be given. Let $A$ be an $m\times n$ simple matrix with no $q\cdot (\1_2\0_1)$. Assume $m\ge 6$ and (\ref{a2a3a4110}).  Then
\begin{E}\binom{2}{2}\binom{m-2}1a_2+
\binom{3}{2}\binom{m-3}1a_3
+\binom{4}{2}\binom{m-4}1a_4\le \binom{m}{3}3(q-1).\label{pigeonholeq110}\end{E}
There exists positive constants $c_1,c_2$ so that
\begin{E}\binom{m}{2}-c_1m\le a_2\le \binom{m}{2},\label{a2ineqq110}\end{E}
\begin{E}a_4\le c_2m.\label{a4ineq110}\end{E}\end{lemma}
\proof We note that a column of column sum $k$ has $\binom{k}{2}\binom{m-k}1$ configurations $\1_2\0_1$ and note that $\binom{k}{2}\binom{m-k}{1}\ge
\binom{4}{2}\binom{m-4}{1}$ for $4\le k\le m-1$.  Counting the configurations $\1_2\0_1$ and using the pigeonhole argument yields (\ref{pigeonholeq110})

For $m\ge 6$ we have $\binom{3}{2}\binom{m-3}1<
\binom{4}{2}\binom{m-4}1$.
Hence
$$(m-2)a_2+3(m-3)(a_3+a_4)\le \binom{m}{3}3(q-1)$$
From (\ref{a2a3a4110}), we have $a_3+a_4\ge \frac{q+1}{3}\binom{m}{2}-a_2$.
We substitute and obtain
$$(m-3)\binom{m}{2}(q+1)-\binom{m}{3}3(q-1)\le 
\biggl(3(m-3)-(m-2)\biggr)a_2$$
which simplifies as 
$$\binom{m}{2}(2m-q-5)\le (2m-7)a_2$$
from which we deduce that there is a constant $c_1$ (will depend on $q$) so that first half of (\ref{a2ineqq110}) holds. The second half of (\ref{a2ineqq110}) follows from the fact that no column of sum 2 is repeated.

In a similar way we have 
$$(m-2)a_2+3(m-3)\biggl(\frac{q+1}{3}\binom{m}{2}-a_2-a_4\biggr)+6(m-4)a_4
\le\binom{m}{3}3(q-1)$$
and when we substitute the upper bound of (\ref{a2ineqq110}), we  deduce that there is a constant $c_2$
(will depend on $q$) so that (\ref{a4ineq110}) holds.\qed

Partition $A$ into three parts: $A_{2}$ consists of the columns of column sum 2, $A_{3}$ is the  columns of column sum 3 and $A_{4}$ is the columns of column sum greater or equal than 4. We will refer to ${\cal A}_2$, ${\cal A}_3$ using the notations of (\ref{salpha}) and (\ref{calA}). Note that ${\cal A}_3$ is a multiset and ${\cal A}_2$ is a set given that there are no repeated columns of sum 2.  Considering the columns of column sum 2, we adapt $\epsilon(i)$ of (\ref{e1e0}).  Note that for convenience we represent every pair $\{i,j\}$ by $ij$ and so $ij\equiv ji$. We are not interested in ordered pairs in this context. Define
$$\epsilon(ij)=\left\{ \begin{array}{ll} 1 & \hbox{ if } \{i,j\} \notin {\cal A}_{2} \\  0 &\hbox{ if } \{i,j\} \in {\cal A}_{2} \end{array} \right.\quad ,\qquad
E=\{ij\,:\,\epsilon(ij)=1\}
\quad .$$
Thus 
\begin{E}a_2=\binom{m}{2}-\sum_{ij}\epsilon(ij)=\binom{m}{2}-|E|.\label{a2eq110}\end{E}
We deduce from (\ref{a2ineqq110}) that $|E|\le c_1m$. 

We adapt the definitions of the degrees $d_0$, $d_1$ of Section~\ref{sectionq10} by using a hypergraph degree definitions applied to the multiset ${\cal A}_3=\{B_1,B_2,\ldots\}$.
Define 
$$d(ij)=|\{s\,:\,B_s\in{\cal A}_3\hbox{ and }i,j\in B_s\}|$$ 
Then
\begin{E}3a_3=\sum_{\{i,j\}\subseteq [m]}d(ij).\label{a3eq}\end{E}
Let 
$${\cal U}(pt)=\{r: \{p,t,r\} \in {\cal A}_{3}\}.$$
 Since $m>q+2$ and we are avoiding $q\cdot (\1_2\0_1)$ in $A_{3}$ then $|{\cal U}(pt)|<q$. Also let 
$${\cal T}(r)=\{pt: \{p,t,r\} \in {\cal A}_{3}\}.$$
 Since for every  $x \in [m]$ with $x\neq r$, $|{\cal U}(rx)|<q$ we have $|{\cal T}(r)|< \frac{(m-1)q}{2}$. Note that ${\cal U}(pt)$ and ${\cal T}(r)$ are the generalizations of $X$ (found after (\ref{dupperbdwitht})) given in the proof of Theorem~\ref{q10}.

\vskip 5pt
\begin{lemma} We have 
\begin{E}d(ij) \leq (q-2)+\epsilon(ij). \label{dijineq}\end{E}\label{dij} \end{lemma}

\proof Since $m>q+2\ge|{\cal U}(ij)|+2$, for every pair $ij$, we can find row $k\neq i,j$ so that $k \notin {\cal U}(ij)$. Now the number of submatrices
  \begin{E} \begin{matrix} i \\j \\k\\ \end{matrix} \begin{bmatrix} 1 \\1 \\0 \end{bmatrix} \label{submatrix110} \end{E} 
in $A_{3}$ is $d(ij)$ (since $d(ij)$ is the number of triples $i,j,l$ corresponding to columns in $A_3$ and each such column yields the submatrix since $k\ne{\cal U}(ij)$) and  the number of submatrices ($\ref{submatrix110}$) in $A_{2}$ is $1-\epsilon(ij)$. Thus $$d(ij)+1-\epsilon(ij)\leq q-1$$ and hence  (\ref{dijineq}) holds.\qed

Let 
$$Y=\{ij\,:\,d(ij)=q-2\hbox{ and }\epsilon(ij)=0\}$$

\begin{lemma}\label{Y110}There exists a constant $c_3$ so that 
\begin{E}|Y|\ge \binom{m}{2}-c_3m\label{yineq}\end{E}\end{lemma}
\proof
We  partition  the $\binom{m}{2}$ pairs $ij$ into 3 parts: $Y$, $E$ and the rest. By Lemma~\ref{dij}, for each $ij\in E$ we have $d(ij)\le (q-2)+1$. Note that for $ij\notin Y\cup E$, we have $\epsilon(ij)=0$ and so $d(ij)\le (q-2)-1$ else $ij\in Y$. Thus from (\ref{a3eq}) 
$$3a_3=\sum_{ij}d(ij)\le\left((q-2)|Y|
+((q-2)+1)|E|
+((q-2)-1)\biggl(\binom{m}{2}-|Y|-|E|\biggr)\right)$$
Hence 
\begin{E}
a_3\le\frac{1}{3}\left((q-2)\binom{m}{2}+|E|-\binom{m}{2}+|Y|+|E|\right)
\label{a3ineq110extend}\end{E}
Substituting estimates of $a_2$, $a_3$, $a_4$ from (\ref{a2eq110}), (\ref{a3ineq110extend}), (\ref{a4ineq110}) into (\ref{a2a3a4110}), we have
$$\binom{m}{2}-|E|+  \frac{1}{3}\left((q-2)\binom{m}{2}+2|E|-\binom{m}{2}+|Y|\right)
+c_2m
>\frac{q+1}{3}\binom{m}{2}$$
We deduce $-\frac{1}{3}|E|+\frac{1}{3}|Y| +c_2m>\frac{1}{3}\binom{m}{2}$ and so there exists a constant $c_3=3c_2$ so that (\ref{yineq}) holds.\qed

Form a graph $G$ of $m$ vertices corresponding to the rows of $A$ and 
with edges $(i,j)$ if and only if $ij\in Y$. Thus by Lemma~\ref{Y110}, the number of edges
of $G$ is at least $\binom{m}{2}-c_3m$. By Tur\'an's Theorem \cite{T},  a graph with more than $\frac{m^{2}}{2}-\frac{m^{2}}{2(k-1)}$ edges has a clique of $k$ vertices. Thus $G$ has large cliques.   Let $c_4$ be a constant chosen so that for any choices of $i,j,k$ the following three inequalities hold.
$$\binom{c_4\sqrt{m}/2}{2}>\frac{m-1}{2}q \bigl(> 
|{\cal T}(k)|\bigr),\quad \frac{c_4{\sqrt{m}}}{2} >q\bigl(>|{\cal U}(ij)|\bigr),$$
\begin{E} \binom{c_4\sqrt{m}}{2}>\frac{m-1}{2}q+3m
\bigl(\ge 
|{\cal T}(k)|+|{\cal U}(ij)| \bigr)\label{clique110}\end{E}
By Tur\'an's argument, there exists an $M$ so that for $m\ge M$, we can find a clique of  $c_4\sqrt{m}$ vertices in $G$.   Let the vertices in this clique be denoted $B$. Thus for $i,j\in B$ we have $d(ij)=q-2$ and $\epsilon(ij)=0$. Let $A_4(B)$ be the submatrix of $A_4$ of the rows indexed by $B$.

\begin{lemma}Assume $m>M$. Then $A_4(B)$ has no configuration $\1_2\0_1$.\label{no110}\end{lemma}
\proof 
Consider a column $\alpha$ of $A_4$. We consider two cases based on whether there are more 1's or more 0's in the rows $B$. Assume $\alpha$ has at least $\frac{c_4{\sqrt {m}}}{2}$ 1's in rows of $B$. Assume $\alpha$ has  a 0 in row $k\in B$. Then by the first inequality (\ref{clique110}), there is a pair $ij\notin {\cal T}(k)$ with $i,j\in B$. Thus there are $q-2$ columns of column sum 3 with the submatrix (\ref{submatrix110}) using $d(ij)=q-2$ and 1 column of column sum 2 with the submatrix (\ref{submatrix110}) using $\epsilon(ij)=0$ and column $\alpha$ has 1 further submatrix (\ref{submatrix110}) which creates the configuration $q\cdot (\1_2\0_1)$, a contradiction. So $\alpha$ has no configuration $\1_2\0_1$.

 Assume $\alpha$ of $A_4$ that has at least $\frac{c_4{\sqrt {m}}}{2}$ 0's in the rows of $B$. Assume $\alpha$ has  1's  in rows $i,j\in B$.  Then there is a row $k\in B$ where $\alpha$ has a 0 in row $k$ and $k\notin {\cal U}(ij)$ by the second inequality of (\ref{clique110}).  For that choice of $k$ and using $d(ij)=q-2$,   there are $q-2$ columns of column sum 3 with the submatrix (\ref{submatrix110}). There is  one column of column sum 2 with the submatrix (\ref{submatrix110}) using $\epsilon(ij)=0$ and the $\alpha$ has one further submatrix (\ref{submatrix110}) which creates the configuration $q\cdot (\1_2\0_1)$, a contradiction. Thus $\alpha$ has no configuration $\1_2\0_1$. \qed

\begin{lemma} Assume $m>M$. Then the inequality  (\ref{q110bd}) holds.\label{onlybd110}\end{lemma}

\proof  We obtain a contradiction from assuming $m>M$ and (\ref{a2a3a4110}) and thus establish (\ref{q110bd}). Our proof considers the $a_4$ columns of $A_4$ (which are the columns of column sum at least 4 and at most $m-1$).

From Lemma~\ref{no110}, each column in $A_4$ either has at most one 1 or has no 0's in the rows of $B$. Let  $A_4^0$ be those columns of $A_4$  with  at most one 1 in the rows of $B$ and hence at least three 1's in the rows $[m]\backslash B$. Let $a_4^0$ be the number of columns in $A_4^0$. Let $A_4^1$ be those columns of $A_4$  with  no 0's in the rows of $B$ and hence at least one 0 in the rows $[m]\backslash B$. 
Let $a_4^1$ be the number of columns of $A_4^1$. We have $a_4^0+a_4^1=a_4$.

For a pair $ij$ with $i,j\in [m]\backslash B$, let $t(ij)$ 
count the number of columns of $A_4^0$ with 1's in both rows $i$ and $j$.   Each column with at most one 1 in $B$ has at least three 1's in $[m]\backslash B$ and hence
1's in at least $\binom{3}{2}=3$ pairs $ij$ with $i,j\in [m]\backslash B$.
We have verified that 
\begin{E}\sum_{ij\,:\,i,j\in [m]\backslash B}t(ij)\ge 3a_4^0.\label{sumtij0}\end{E}

We must work harder to get an analog of (\ref{sumtij0}) for $A_4^1$. Assume $a_4^1>0$.
For a pair $ij$ with $i,j\in B$ and $k\in [m]\backslash B$ with $ij\notin{\cal T}(k)$, let
$t(ij,k)$ denote the number of submatrices (\ref{submatrix110}) in  $A_4^1$. When $ij\in{\cal T}(k)$,  set $t(ij,k)=0$.
For a pair $ij$ with $i,j\in B$, let 
\begin{E}t(ij)=\max_{k\in [m]\backslash B}t(ij,k)\label{tijdefn}\end{E}

 Each column $\alpha$ in $A_4^1$ has at least one row, say $l\in [m]\backslash B$ with a 0.  For  column $\alpha$, we know
$|{\cal T}(l)|<\frac{m-1}{2}q$ and at the same time there are
$\binom{c_4{\sqrt m}}{2}$ pairs $ij$ with $i,j\in B$ and so there are at least
$\binom{c_4{\sqrt m}}{2}-\frac{m-1}{2}q$ pairs $ij$ with $i,j\in B$ with 
$ij\notin {\cal T}(l)$. Thus by the third inequality of (\ref{clique110}), column $\alpha$ contributes at least $3m$ to the sum
$\sum_{ij\notin{\cal T}(l)}t(ij,l)$ and so
$$\sum_{l\in [m]\backslash B}\quad\sum_{ij\,:\,i\,j\in B}t(ij,l)>3ma_4^1$$
Thus by (\ref{tijdefn}), 
$$(m-|B|)\cdot\sum_{ij\,:\,i,j\in B}t(ij)> \sum_{l\in [m]\backslash B}\quad\sum_{ij\,:\,i\,j\in B}t(ij,l)$$
and so we deduce that
\begin{E}\sum_{ij\,:\,i,j\in  B}t(ij)>\frac{3ma_4^1}{m-|B|}> 3a_4^1.\label{sumtij1}\end{E}

For a pair $ij$ with $i\in B$ and $j\in[m]\backslash B$ or vice versa, let $t(ij)=0$.
We add (\ref{sumtij0}) and (\ref{sumtij1}) together to get
\begin{E}\sum_{ij}t(ij)\ge 3a_4,\label{sumtij}\end{E}
with strict inequality if $a_4^1>0$.

We are able to extend Lemma~\ref{dij} and establish
\begin{E}d(ij)\le q-2+\epsilon(ij)-t(ij)\label{dupperbdwithtij}\end{E} 
By Lemma~\ref{dij}, we need only consider $ij$ with $t(ij)>0$. Given the definition of $t(ij)$, we need only consider the two cases:
$i,j\in [m]\backslash B$ or $i,j\in B$. 

In the former case we note that each of the $t(ij)$ columns of $A_4^0$ with 1's in both rows $i$ and $j$ have at most one 1 in rows of $B$. With $|B|>2q$, (by the second inequality of (\ref{clique110})) we deduce that $t(ij)<q$ else we will find the configuration $q\cdot (\1_2\0_1)$ in $A_4^0$ in the rows $i,j$ and a row of $B$. Now in these $t(ij)$ columns of $A_4^0$, at least $|B|-q+1$ rows of $B$ are all 0's.  Again using the second inequality of (\ref{clique110}) that $|B|>2q$, we can find some $k\in B$ with $k\notin {\cal U}(ij)$ and all the $t(ij)$ columns have 0's in row $k$. 
Now there are 
$d(ij)$ submatrices (\ref{submatrix110}) in columns of sum 3, 
$(1-\epsilon(ij))$ submatrices (\ref{submatrix110}) in columns of sum 2, and
$t(ij)$ submatrices (\ref{submatrix110}) in columns of sum 4 or more. The total is at most $q-1$ since otherwise we would have the configuration $q\cdot (\1_2\0_1)$ and this yields $d(ij)+(1-\epsilon(ij))+t(ij)\le q-1$. This is (\ref{dupperbdwithtij}). 

In the  latter case with $i,j\in B$, we select $k$ so that $t(ij,k)= t(ij)$. Thus $k\notin {\cal U}(ij)$ and also there are at least $t(ij)$ submatrices (\ref{submatrix110}) in columns of $A_4^1$.  Thus we can now follow the same argument as in the former case to establish (\ref{dupperbdwithtij}).

Now using (\ref{dupperbdwithtij}) and (\ref{sumtij}),
\begin{E} 3a_3=\sum_{ij} d(ij) \leq \sum_{  i,j} \biggl(q-2+\epsilon(ij) -t(ij)\biggr)=(q-2)\binom{m}{2}+|E|-3a_4. \label{sumdij} \end{E}
 Substituting (\ref{a2eq110}), (\ref{sumdij}),
and (\ref{a4ineq110}) in (\ref{a2a3a4110})  we obtain
$$\binom{m}{2}-|E|+ \frac{1}{3}\left((q-2)\binom{m}{2}+|E|-3a_4 \right) +a_4>\frac{q+1}{3}\binom{m}{2}.$$
Simplifying and rearranging,
$$-\frac{2}{3}|E|>0$$
which is a contradiction (even for $|E|=0$) and this establishes (\ref{q110bd}). \qed

\vskip 5pt
\noindent{\bf Proof of Proposition~\ref{nonsimpleq110}}:
Lemma~\ref{onlybd110} establishes most of Proposition~\ref{nonsimpleq110} but we are also interested in cases when the bound is achieved.
Assume $m>M$ and $m\equiv 1,3(\hbox{mod }6)$.
We now consider an $m$-rowed simple matrix $A$ which has no configuration $q\cdot (\1_2\0_1)$ and with
 $\binom{m}{0}+\binom{m}{1}+\frac{q+1}{3}\binom{m}{2}+\binom{m}{m}$ columns. One repeats  the previous lemmas and arguments replacing the inequality (\ref{a2a3a4110}) with the equation
\begin{E}a_2+a_3+a_4=\frac{q+1}{3}\binom{m}{2}.\label{a2a3a4eq110}\end{E}
 
We wish to show $a_2=\binom{m}{2}$, $a_4=0$, $a_3=\frac{q-2}{3}\binom{m}{2}$. 
Now Lemma~\ref{pigeon110} holds with (\ref{a2a3a4110}) as an equality.
We deduce the same bounds for ${\cal U}(rx)$ and ${\cal T}(r)$. Lemma~\ref{dij} still holds since the final contradiction does not require the strict inequality of (\ref{a2a3a4110}) merely the equality of
(\ref{a2a3a4eq110}). Lemma~\ref{Y110} holds and we can choose $B$ as large as possible but at least satisfying the three inequalities (\ref{clique110}).
Lemma~\ref{no110} continues to hold.

We use (\ref{a2a3a4eq110}) and following the  argument of Lemma~\ref{onlybd110}, we deduce that $E=\emptyset$ and so $a_2=\binom{m}{2}$. Also we deduce that 
$$\sum_{ij}t(ij)=3a_4$$
and as a result of the strict inequality in (\ref{sumtij1}), we can deduce that $a_4^1=0$. 

Assume $a_4=a_4^0>0$ and consider $\alpha$ in $A_4$ with column sum 4 and with 1's in rows $i,j,k,l$ where $i\in B$ and 
$j,k,l\in\{1,2,\ldots ,m\}\backslash B$. Choose $r\in B\backslash i$ then $\alpha$ has 1's in rows $i,j$ and 0's in row $r$. Using $E=\emptyset$, we deduce that for this particular $i,j$ we have 
$d(ij)\le (q-2)-1$. This yields a slight variant of (\ref{sumdij}): 
$$3a_3=\sum_{ij}d(ij) \leq 
\sum_{  ij}{((q-2)+\epsilon(ij)-t(ij))}-1.$$
The extra `-1' is sufficient to obtain a contradiction when we substitute for $a_2,a_3,a_4$ in (\ref{a2a3a4eq110}). We then deduce $a_4=0$.

With $a_4=0$ and $a_2=2\binom{m}{2}$, we deduce
$a_3=\frac{q-3}{3}\binom{m}{2}$ using (\ref{a2a3a4eq110}). Given that $\epsilon(ij)=0$ for all $ij$ and using Lemma~\ref{dij}, we deduce 
$d(ij)=q-2$ for all pairs $ij$ and so $B=\{1,2,\ldots ,m\}$. From this we can readily conclude that the columns of column sum 3 correspond to a 2-design $S_{q-2}(2,3,m)$ and $A$ has no further columns. \qed

\section{Exact Bound for $q\cdot(\1_2\0_2)$}\label{sectionq1100}

We generalize our proof of Proposition~\ref{nonsimpleq110} given in Section~\ref{sectionq110} to prove Proposition~\ref{nonsimpleq1100}. Again we do not explicitly calculate the smallest possible constant $M$  but we note that we can take $M$ to be $O(q^3)$.

Let $A$ be a $m\times n$ matrix with no $q\cdot(\1_2\0_2)$. Assume that there are no repeated columns of sums $0,1,2,m-2,m-1,m$. We will assume $n > 2+2m+\binom{m}{2}\frac{q+3}{3}$. Let $a_{i}$ denote the number of columns with either exactly $i$ 1's or $i$ 0's for $i=0,1,2,3$ and let $a_{4}$ be the number of remaining columns. We may assume $a_0=2$ and $a_1=2m$ since all columns of column sum $0,1,m-1$ or $m$ do not contain the configuration $\1_2\0_2$.  Thus 
\begin{E} a_{2}+a_{3}+a_{4}>\binom{m}{2}\frac{q+3}{3}. \label{a2a3a41100} \end{E}

\begin{lemma} Assume $A$ is an $m\times n$ simple matrix with no configuration $q\cdot(\1_2\0_2)$ and (\ref{a2a3a41100}) holds. Then there exists an $m_0$ so that for $m>m_0$,
\begin{E}\binom{2}{2}\binom{m-2}{2}a_{2}+\binom{3}{2}\binom{m-3}{2}a_{3}+\binom{4}{2}\binom{m-4}{2}a_{4} \leq 6\binom{m}{4}(q-1) .\label{inequality1} \end{E}
Also there exist constants $c_1,c_2$ so that
\begin{E} 2\binom{m}{2}-c_1m\le a_2\le 2\binom{m}{2}\label{a2ineq1100}\end{E}
\begin{E}a_{4} \leq c_{2}m\label{a4ineq1100}\end{E} 
\label{pigeon1100}\end{lemma}
\proof
A column in $A$ of column sum $k$ has $\binom{k}{2}\binom{m-k}{2}$ configurations $\1_2\0_2$. Note that $\binom{k}{2}\binom{m-k}{2}\ge
\binom{4}{2}\binom{m-4}{2}$ for $4\le k\le m-4$. By the pigeonhole principle, there are at most $6(q-1)\binom{m}{2}$ configurations $\1_2\0_2$ in $A$. We obtain (\ref{inequality1}).
There exist an $m_0$, such that for $m>m_0$, $\binom{3}{2}\binom{m-3}{2} <\binom{4}{2}\binom{m-4}{2}$. Substituting in $(\ref{inequality1})$,
$$ \binom{m-2}{2}a_{2}+3\binom{m-3}{2}(a_{3}+a_{4}) \leq 6(q-1)\binom{m}{4} $$
which yields using $a_{3}+a_{4} > \binom{m}{2}\frac{q+3}{3}-a_2$ and rearranging
$$ \binom{m-2}{2}a_{2}+3\binom{m-3}{2}\left [\binom{m}{2}\frac{q+3}{3}-a_{2}\right]< 6(q-1)\binom{m}{4}.$$
Therefore,
\begin{E} \binom{m-3}{2}\binom{m}{2}(q+3)-6(q-1)\binom{m}{4} < \left(3\binom{m-3}{2}-\binom{m-2}{2}\right)a_{2}. \label{a2} \end{E}
The leading term on the righthand side is exactly $m^4$ while the leading coefficient of $a_2$ on the lefthand side is exactly $m^2$. Thus (\ref{a2}) implies that there exists some constant $c_1$ so that the lower bound of (\ref{a2ineq1100}) holds. The upper bound of
(\ref{a2ineq1100}) follows from the fact no column of sum 2 or $m-2$ is repeated. 

We can also bound $a_{4}$. From (\ref{a2a3a41100}), we have $a_3>\frac{q+3}{3}\binom{m}{2}-a_2-a_4$.  Using $(\ref{inequality1})$ we have
$$ \binom{m-2}{2}a_{2}+3\binom{m-3}{2} \left [\frac{q+3}{3}\binom{m}{2}-a_{2}-a_{4} \right ]+6\binom{m-4}{2}a_{4} \leq 6(q-1)\binom{m}{4} .$$
Then
$$
\left[6\binom{m-4}{2}-3\binom{m-3}{2}\right]a_{4} $$
$$\leq 6(q-1)\binom{m}{4}-
(q+3)\binom{m-3}{2}\binom{m}{2}
+a_{2}\left[3\binom{m-3}{2}-\binom{m-2}{2}\right] .$$
Substituting $a_{2} \leq 2 \binom{m}{2}$ and rearranging we have
\begin{E} \left[6\binom{m-4}{2}-3\binom{m-3}{2}\right]a_{4} \leq \frac{m(m-1)(m-3)}{4}(2q-6) \label{a4} \end{E}
Then $(\ref{a4})$ implies that there exist some constant $c_2$ so that (\ref{a4ineq1100}) holds.\qed

We could have produced the bound $a_4\le (2q-6)\frac{m}{6}+c_2'$ for some constant $c_2'$, but this is of little help. Now we form analogs of the degrees
 $d_0$, $d_1$ of Section~\ref{sectionq10} by defining $A_3$ as the submatrix of $A$ of the columns of column sum 3 and defining $A_{m-3}$ as the submatrix of $A$ of the columns of column sum $m-3$. We refer to the mutisets ${\cal A}_3=\{B_1,B_2,\ldots\}$, ${\cal A}_{m-3}=\{C_1,C_2,\ldots\}$ using the notations of (\ref{salpha}) and (\ref{calA}). Define
$$d_1(ij)=|\{s\,:\,B_s\in {\cal A}_3\hbox{ and }i,j\in B_s\}|,\quad 
d_0(ij)=|\{s\,:\,C_s\in {\cal A}_{m-3}\hbox{ and }i,j\notin C_s\}|$$  
Recalling $a_3=|{\cal A}_3|+|{\cal A}_{m-3}|$, we note
\begin{E} 3a_3= \sum_{ \{i,j\}\subset [m]}\bigl(d_0(ij)+d_1(ij)\bigr) \label{a3eq1100}\end{E}

Define $e_{ij}$ to be the $m$-rowed column with 1 in rows $i$ and $j$ and 0's elsewhere, and let $e_{ij}^{c}$ be the (0,1)-complement of  $e_{ij}$. These are the possible columns of column sum 2 or $m-2$. Define
$$ E_1=\{ij: { \{i,j\}\subset [m]} \hbox{ and } e_{ij} \hbox{ is not in } A\} $$
$$   E_0=\{ij: { \{i,j\}\subset [m]} \hbox{ and } e_{ij}^{c} \hbox{ is not in } A\} $$
For convenience of counting define
\begin{E}\epsilon(ij)=\left\{\begin{array}{ll}0&\hbox{ if }ij\notin E_1\cup E_0\\
1&\hbox{ if }ij\in E_1\backslash E_0\hbox{ or }ij\in E_0\backslash E_1\\ 
2&\hbox{ if }ij\in E_1\cap E_0\\\end{array}\right.\quad .\label{e1100}\end{E}

Thus
\begin{E}a_{2}=2\binom{m}{2}- \sum_{ i,j\subset [m]}{\epsilon(ij)}= 2\binom{m}{2}-(|E_1|+|E_0|),\label{a2eq1100}\end{E}
and given (\ref{a2ineq1100}) we have $|E_1|+|E_0| \le c_1m$

We note for a quadruple of rows $p,t,r,s$ that are at most $2q-2$
\begin{E} \hbox{ submatrices } 
\begin{matrix} p\\ t\\ r\\ s\\ \end{matrix} 
\begin{bmatrix} 1 \\ 1 \\ 0 \\ 0\\ \end{bmatrix} 
\hbox{ or submatrices } 
\begin{matrix} p\\ t\\ r\\ s\\ \end{matrix} 
\begin{bmatrix} 0 \\ 0 \\ 1 \\ 1\\ \end{bmatrix} \label{submatrices} \end{E}
else $A$ has the configuration $q\cdot (\1_2\0_2)$. For disjoint pairs $pt$ and $rs$ (i.e. $\{p,t\}\cap \{r,s\}=\emptyset$) 
we say  {\it pair $pt$ has triple overlapping $rs$} if and only if at least one of submatrices 
$$\left.\begin{matrix} p \\ t\\ r \\ s \end{matrix}\right. \begin{bmatrix} 1 \\ 1 \\ 1 \\ 0 \end{bmatrix}  \hbox{ or }
\left.\begin{matrix} p \\ t\\ r \\ s \end{matrix}\right. \begin{bmatrix} 1 \\ 1 \\ 0 \\ 1 \end{bmatrix} $$ 
appears in columns of column sum 3 or at least one of submatrices 
$$\left.\begin{matrix} p \\ t\\ r \\ s \end{matrix}\right.\begin{bmatrix} 0 \\ 0 \\ 0 \\ 1 \end{bmatrix} \hbox{ or } 
\left.\begin{matrix} p \\ t\\ r \\ s \end{matrix}\right. \begin{bmatrix} 0 \\ 0 \\ 1 \\ 0 \end{bmatrix} $$ 
appears in columns of column sum $m-3$. This definition is not symmetric in the pair $pt,rs$. Note that columns of three 1's that have 1's on rows $p,t$ yet no 1's on rows $r,s$ or vice versa have 1's on rows $r,s$  yet no 1's on rows $p,t$ contribute to (\ref{submatrices}). Similarly for columns with three 0's. Let
$$ {\cal U}(pt)=\{ij: { \{i,j\}\subset [m]}, \hbox{pair } pt \hbox{ has triple overlapping } ij \}, $$
$$ {\cal T}(pt)=\{ij: { \{i,j\}\subset [m]}, \hbox{pair } ij \hbox{ has triple overlapping } pt \} .$$
Given $m>q+2$, we cannot have the submatrix 
$q\cdot(\1_2\0_0)$ in rows $p,t$
in columns of column sum 3 else we would have the configuration $q\cdot(\1_2\0_2)$ (and so there are at most $q-1$ columns of column sum 3 with 1's in rows $p,t$).  Similarly, we  cannot have the submatrix
$q\cdot(\1_0\0_2)$ in rows $p,t$ 
in columns of column sum $m-3$. To bound ${\cal U}(pt)$, we note that 
$\binom{m-2-(q-1)}{2}$ counts the number of pairs $ij$ disjoint from $pt$ that avoids $q-1$ further rows. Thus the number of pairs $ij$ where $pt$ overlaps $ij$ using a column of column sum 3 is at most $\binom{m-2}{2}-\binom{m-2-(q-1)}{2}$. Similarly, the number of pairs $ij$ where $pt$ overlaps $ij$ using a column of column sum $m-3$ is at most $\binom{m-2}{2}-\binom{m-2-(q-1)}{2}$. Thus there exists a constant $c_3$ depending only on $q$ so that
\begin{E} |{\cal U}(pt)| \leq 2\bigl(\binom{m-2}{2}-\binom{m-2-(q-1)}{2}\bigr)\le c_3m .\label{u} \end{E}
Given $m>q+2$ and a fixed choice $x$ different from $p,t$, we note that the columns of column sum 3 cannot have 
the submatrix $q\cdot(\1_2\0_0)$ in rows $p,x$ nor the submatrix $q\cdot(\1_2\0_0)$ in rows $t,x$ since either would produce the configuration
$q\cdot (\1_2\0_2)$.
Thus  for a fixed $x\ne p,t$ (of which there are $m-2$ choices), there are at most $2(q-1)$ choices for $j$ such that pair $xj$ has triple overlapping $pt$ in columns of column sum 3. A similar argument applies to the columns of column sum $m-3$. Thus there exists a constant $c_4=2(q-1)$  so that
\begin{E} |{\cal T}(pt)| \leq 2\left(\frac{(m-2)2(q-1)}{2}\right)\leq c_4m \label{t} \end{E}

\vskip 5pt
\begin{lemma} There exists a constant $m_1\ge q+4$ so that for $m>m_1$, we have for all ${ \{i,j\}\subset [m]}$ that $d_0(ij)+d_1(ij) \leq q-3+\epsilon(ij)$ .\label{dupperbd1100}\end{lemma}
\proof
Assume the contrary that $pt$ is an index with $d_0(pt)+d_1(pt) \geq q-3+\epsilon(pt)+1$. Let $\{r,s\}\subset [m]\backslash \{p,t\}$ and $rs \notin {\cal U}(pt)\cup {\cal T}(pt)$ (by $(\ref{u})$ and $(\ref{t})$ there are $\binom{m}{2}-c_3m-c_4m$ choices for $rs$). There are at most $2q-2$
submatrices as in (\ref{submatrices}) 
contained in $A$ else $A$ has the configuration $q\cdot (\1_2\0_2)$. There are $4-\epsilon(pt)-\epsilon(rs)$ submatrices $(\ref{submatrices})$ contained in columns of column sum 2, $m-2$ and since $rs \notin {\cal U}(pt)\cup {\cal T}(pt)$ there are $(d_1(pt)+d_0(rs))+(d_0(pt)+d_1(rs))$ submatrices $(\ref{submatrices})$ in columns of column sum 3, $m-3$. Thus
$$ (d_1(pt)+d_0(rs))+(d_0(pt)+d_1(rs))+4-\epsilon(pt)-\epsilon(rs) \leq 2(q-1) $$
Substituting $d_0(pt)+d_1(pt) \geq q-3+\epsilon(pt)+1$ and rearranging yields
\begin{E} d_0(rs)+d_1(rs)  \le (q-3)-1+\epsilon(rs). \label{rsinequality} \end{E}

We wish to bound $a_3$ using (\ref{a3eq1100}). We split all pairs $ij$ into three sets: those with $\{i,j\}\cap \{p,t\}=\emptyset$ and $ij\notin {\cal U}(pt)\cup {\cal T}(pt)$, those with
$ij\in {\cal U}(pt)\cup {\cal T}(pt)$ (which forces $\{i,j\}\cap\{p,t\}=\emptyset$) and those with $\{i,j\}\cap \{p,t\}\ne\emptyset$. In the first case, we use (\ref{rsinequality}). 

$$\sum_{ \begin{subarray}{c}\{i,j\}\subset [m]\\
ij\notin {\cal U}(pt)\cup {\cal T}(pt) \\
 \{i,j\}\cap \{p,t\}=\emptyset\\ \end{subarray}}  d_0(ij)+d_1(ij)
\le 
\sum_{ \begin{subarray}{c}\{i,j\}\subset [m]\\
ij\notin {\cal U}(pt)\cup {\cal T}(pt) \\
 \{i,j\}\cap \{p,t\}=\emptyset\\ \end{subarray}}
 (q-3)-1+\epsilon(ij)$$

In the latter cases, note that $d_0(ij)\le q-1$ and $d_1(ij)\le q-1$ else, since $m\ge q+4$, we would find a copy of $q\cdot (\1_2\0_2)$.

$$\sum_{\begin{subarray}{c}\{i,j\}\subset [m]\\
ij\in {\cal U}(pt)\cup {\cal T}(pt)\\ \end{subarray}} d_0(ij)+d_1(ij)
\le (c_3+c_4)m\cdot 2(q-1)$$

$$\sum_{\begin{subarray}{c}\{i,j\}\subset [m] \\ \{i,j\}\cap \{p,t\}\neq \emptyset\end{subarray}}    
 d_0(ij)+d_1(ij)
\le
2(m-2)\cdot 2(q-1)$$
Let $c_5$ be a constant chosen so that $c_5>2(c_3+c_4+2)(q-1)$. Combining yields
$$\sum_{ij}(d_0(ij)+d_1(ij))\le \sum_{\begin{subarray}{c} \{i,j\}\subset [m]\\
ij\notin {\cal U}(pt)\cup {\cal T}(pt) \\
 \{i,j\}\cap \{p,t\}=\emptyset\\\end{subarray} }
 \biggl((q-3)-1+\epsilon(ij)\biggr)+c_5m.$$

Now using (\ref{a2a3a41100}) and substituting for $a_2$ using (\ref{a2eq1100}) and substituting for $a_3$ using (\ref{a3eq1100}) and the above inequality with the estimate that there are at most
$\binom{m-2}{2}$ choices for pairs $ij$ with $\{i,j\}\cap\{p,t\}=\emptyset$
$ij\notin{\cal U}(ij)\cup{\cal T}(ij)$ and substituting for $a_4$ using (\ref{a4ineq1100}):

$$2\binom{m}{2}-(|E_0|+|E_1|) +
 \frac{(q-3)-1}{3}\binom{m-2}{2}+\frac{1}{3}(|E_0|+|E_0|)+\frac{c_5}{3}m+c_2m>\binom{m}{2}\frac{q+3}{3}
$$
The coefficient of $m^{2}$ on the left side of the above inequality is only $\frac{q+2}{6}$ while on the right side is $\frac{q+3}{6}$. Thus there exists a constant $m_1$ so that for $m>m_1$, we have  a contradiction proving the claim.\qed

Let 
$$Y=\{ij\,:\,d_0(ij)+d_1(ij)=q-3\hbox{ and }\epsilon(ij)=0\}$$

\begin{lemma}\label{Y}  There exists a constant $c_6$ so that 
\begin{E}|Y| > \binom{m}{2}-c_{6}m.\label{yineq1100}\end{E}\end{lemma} 
\proof
We partition  the $\binom{m}{2}$ pairs $ij$ into 3 parts: $Y$, $E_0\cup E_1$ and the rest. We note that for $ij\notin Y\cup E_0\cup E_1$, we have 
$\epsilon(ij)=0$ and $d_0(ij)+d_1(ij)\le (q-3)-1$ by Lemma~\ref{dupperbd1100}. Thus from (\ref{a3eq1100}) and using Lemma~\ref{dupperbd1100}
$$a_3=\frac{1}{3}\sum_{ij} d_0(ij)+d_1(ij)\le\frac{1}{3}\left((q-3)|Y|
+((q-3)+2)|E_0\cup E_1|\right.$$
$$\left.+((q-3)-1)\bigl(\binom{m}{2}-|Y|-|E_0\cup E_1|\bigr)\right)$$
Thus
\begin{E}
a_3\le\frac{1}{3}\left((q-3)\binom{m}{2}+3|E_0\cup E_1|-\binom{m}{2}
+|Y|\right)
\label{a3ineq1100}\end{E}
Using (\ref{a2eq1100}),(\ref{a4ineq1100}), (\ref{a3ineq1100}) in (\ref{a2a3a41100}), we have
$$2\binom{m}{2}-(|E_0|+|E_1| )+  \frac{1}{3}\left(
(q-3)\binom{m}{2}+3|E_0\cup E_1|+\biggl(|Y|-\binom{m}{2}\biggr)\right)
+c_2m$$
$$>\frac{q+3}{3}\binom{m}{2}.$$
We deduce, noting that $|E_0|+|E_1| \geq|E_0\cup E_1|$, that 
$\frac{1}{3}\left(|Y|-\binom{m}{2}\right)+c_2m>0$ and so \hfil\break $|Y|>\binom{m}{2}-3c_2m$. Thus (\ref{yineq1100}) holds for $c_6=3c_2$.\qed

\vskip 5pt
Form a graph $G$ whose vertex set is the rows of the matrix $A$ with edges $ij$ for those $ij\in Y$. Thus $G$ has at least $\frac{m^{2}}{2}-c_{6}m$ edges. By Tur\'an's Theorem \cite{T},  a graph with more than $\frac{m^{2}}{2}-\frac{m^{2}}{2(k-1)}$ edges has a clique  of $k$ vertices. Choose  a constant $c_7$ so that for any choices $i,j\in [m]$
$$\binom{c_7{\sqrt m}-2(q-1)}{2}>(c_3+c_4)m\bigl(>|{\cal T}(ij)|+|{\cal U}(ij)|\bigr),$$
$$ \frac{1}{2}\binom{c_7{\sqrt m}}{2}-2m>(c_3+c_4)m\bigl(>|{\cal T}(ij)|+|{\cal U}(ij)|\bigr),$$
\begin{E}\binom{\frac{c_7\sqrt{m}}{2}}{2}>(c_3+c_4)m
\bigl(> |{\cal T}(ij)|+|{\cal U}(ij)|\bigr). \label{clique1100}\end{E}
Then by Tur\'an's Theorem, there exists a $M>m_0,m_1$ ($m_0$ is from Lemma~\ref{pigeon1100} and $m_1$ is from Lemma~\ref{dupperbd1100})so that for $m>M$, graph $G$ has a clique of $c_7{\sqrt{m}}$ vertices.

Let $B$ denote the  set of the rows in this clique. Hence for every $i,j \in B$ we have $ d_1(ij)+d_0(ij)=q-3 \hbox{ and }\epsilon(ij)=0 $. Let $A_4$ denote the columns of $A$ of column sum $4,5,\ldots ,m-5$ or $m-4$.
Let $A_4(B)$ be the submatrix of $A_4$ of the rows indexed by $B$.

\begin{lemma}Assume $m>M$. Then $A_4(B)$ has no configuration $\1_2\0_2$.\label{no1100}\end{lemma}
\proof
Assume there are rows $i,j,k,l\in B$ and  a column $\alpha$ of $A_4$ with 0's in rows $i,j$ and 1's in rows $k,l$. Without loss of generality, we may assume that there are more 1's than 0's in $\alpha$ in the rows of $B$ so that the number of 1's in the rows of $B$ is more than $c_7{\sqrt m}/2$.
Thus by the third inequality in (\ref{clique1100}), we can find a pair $gh$ of rows with $g,h\in B$, so that $\alpha$ has 1's in row $g,h$ and
$gh\notin{\cal T}(ij)\cup{\cal U}(ij)$. We may now argue that for our choice of $i,j,g,h$, we have
$(d_1(ij)+d_0(gh))+(d_1(gh)+d_0(ij))+4-\epsilon(ij)-\epsilon(gh)=2(q-1)$ submatrices
\begin{E}\begin{matrix} i \\ j \\ g \\ h\\ \end{matrix} 
\begin{bmatrix} 1 \\ 1 \\ 0 \\ 0\\ \end{bmatrix}\quad \hbox{ or }\quad
\begin{matrix} i \\ j \\ g \\ h\\ \end{matrix} 
\begin{bmatrix} 0 \\ 0 \\ 1 \\ 1 \end{bmatrix} \label{2submatrices}\end{E}

in $A$ in columns of column sum $2,3,m-3,m-2$. With another such submatrix in $\alpha$ in $A_4$, we have $2(q-1)+1$ such submatrices, for our chosen quadruple $i,j,g,h$ and so $A$ has the configuration $q\cdot(\1_2\0_2)$, a contradiction. \qed

\begin{lemma} Assume $m>M$. Then the inequality  (\ref{q1100bd}) holds.\label{onlybd1100}\end{lemma}

\proof 
Assume $m>M$ and (\ref{a2a3a41100}). Using Lemma~\ref{no1100}, the columns of $A_4$ can be partitioned into two parts: $Z$ the columns that have at most one 1 in the rows $B$ and $J$ the columns that have at most one 0 in the rows of section $B$. 

For each pair $i,j\in[m]\backslash B$, let $t(ij)$ count the sum of the number of columns in $Z$ with 1's in both rows $i,j$ as well as the number of columns in $J$ with 0's in both rows $i,j$. For all other pairs $ij$, let $t(ij)=0$. Given the column sums in $A_4$, every column  in $Z$ has at least three 1's in rows $[m]\backslash B$ and every column in $J$ has at least three 0's in rows $[m]\backslash B$. We have
\begin{E}\sum_{ij}t(ij)\ge 3a_4\label{tsum}\end{E}
Moreover, we find that $t(ij)\le 2(q-1)$: Given a choice for $i,j$, if we have $q$ columns in $Z$ with 1's in rows $i,j$ then there are at most $q$ rows of $B$ containing 1's for these $q$ columns (since each column of $Z$ has at most one 1 in the rows of $B$). But then if we choose two rows of $B$ from the remaining $\ge|B|-q$ rows  in conjunction with $i,j$ then we have a copy of the configuration $q\cdot (\1_2\0_2)$. Similarly, there cannot be $q$ columns of $J$ with 0's on rows $i,j$. We conclude $t(ij)\le 2(q-1)$.

For a given pair $i,j\in [m]\backslash B$, consider  the $t(ij)$ columns contributing to $t(ij)$. By the first inequality in (\ref{clique1100}), 
we can find a pair of rows $gh$ ($g,h \in B$) so that 
$gh\notin{\cal T}(ij)\cup{\cal U}(ij)$ 
and in addition $g,h$ are not chosen from the up to $2(q-1)$ rows of $B$ which are given as follows: the $\le q-1$ rows of $B$ which have 1's in the columns of $Z$ having 1's in both rows $i,j$ and the $\le q-1$ rows of $B$ which have 0's in the columns of $J$ having 0's in both rows $i,j$. Thus if $\alpha$ is a column of $Z$  with 1's in rows $i,j$  then $\alpha$ has 0's in rows $g,h$ and if $\alpha$ is a column of $J$  with 0's in rows $i,j$  then $\alpha$ has 1's in rows $g,h$. 
There will be $4-\epsilon(ij)-\epsilon(gh)$ submatrices as in (\ref{2submatrices}) in the columns of column sum 2 or $m-2$. Neither pair $ij$ has triple overlapping $gh$ nor pair $gh$ has triple overlapping $ij$ and so there will be $(d_1(ij)+d_0(gh))+(d_0(ij)+d_1(gh))$ submatrices as in (\ref{2submatrices}) in the columns of column sum 3 or $m-3$. 
By our choice of $g,h$, a column $\alpha$ in $Z$ with 1's in  rows $i,j$ will have 0's on rows $g,h$.
A column $\beta$ in $J$ with 0's in  rows $i,j$ will have 1's on rows $g,h$.
Thus in $A_4$ we can find $t(ij)$ submatrices as in (\ref{2submatrices}).
In the matrix $A$, an ordered quadruple of rows $i,j,g,h$ has at most $2(q-1)$ submatrices as given in (\ref{2submatrices}) else $A$ would have the configuration
$q\cdot (\1_2\0_2)$.
 Thus
$$ (d_1(ij)+d_0(gh))+(d_0(ij)+d_1(gh))+4-\epsilon(ij)-\epsilon(gh)+t(ij)\leq 2(q-1). $$
Substituting $d_0(gh)+d_1(gh)=q-3$ and $\epsilon(gh)=0$ and rearranging we have
\begin{E} d_1(ij)+d_0(ij) \leq (q-3)+\epsilon(ij)-t(ij).\label{dupperbd1100witht}\end{E}
This inequality is true for other $i,j$ using Lemma~\ref{dupperbd1100} when $t(ij)=0$. Thus
\begin{E}\sum_{ij}(d_0(ij)+d_1(ij)) \leq \sum_{  ij}{(q-3+\epsilon(ij)-t(ij))}\label{sumd0d1}\end{E}
Taking (\ref{a2a3a41100}) with $a_2$ from (\ref{a2eq1100}) and with $a_3$ from
(\ref{a3eq1100}) using (\ref{dupperbd1100witht}) we obtain
$$2\binom{m}{2}
-|E_0|-|E_1|+
\frac{1}{3}\sum_{ij}\bigl(q-3+\epsilon(ij)-t(ij) \bigr)+a_{4}> \frac{q+3}{3}\binom{m}{2}$$
Simplifying and using $\sum_{ij}\epsilon(ij)= |E_0|+|E_1|$ and (\ref{tsum}) we obtain
$$
 -\frac{2}{3}(|E_0|+|E_1|)>0
$$
which is a contradiction (even for $|E_0|+|E_1| =0$). This establishes (\ref{q1100bd}). \qed

\vskip 5pt
\noindent{\bf Proof of Proposition~\ref{nonsimpleq1100}}:
Lemma~\ref{onlybd1100} establishes most of Proposition~\ref{nonsimpleq1100} but we are also interested in cases when the bound is achieved.
Assume $m>M$ and $m\equiv 1,3(\hbox{mod }6)$.
We now consider an $m$-rowed simple matrix $A$ which has no configuration $q\cdot (\1_2\0_2)$ and with
 $\binom{m}{0}+\binom{m}{1}+\frac{q+3}{3}\binom{m}{2}
+\binom{m}{m-1}+\binom{m}{m}$ columns. One repeats  the previous lemmas and arguments replacing the inequality (\ref{a2a3a41100}) with the equation
\begin{E}a_2+a_3+a_4=\frac{q+3}{3}\binom{m}{2}.\label{a2a3a4eq1100}\end{E}

We wish to show $a_2=2\binom{m}{2}$, $a_4=0$, $a_3=\frac{q-3}{3}\binom{m}{2}$ and there exists positive integers $a,b$, $a+b=q-3$ so that for all pairs
$ij$, $d_0(ij)=a$ and $d_1(ij)=b$.
Now Lemma~\ref{pigeon1100} holds with (\ref{a2a3a41100}) as an equality.
We deduce the same bounds for ${\cal U}(pt)$ and ${\cal T}(pt)$. Lemma~\ref{dupperbd1100} still holds since the final contradiction does not require the strict inequality of (\ref{a2a3a41100}) merely the equality of
(\ref{a2a3a4eq1100}). Lemma~\ref{Y} holds and we can choose $B$ as large as possible but at least satisfying the  inequalities (\ref{clique1100}).
Lemma~\ref{no1100} continues to hold. 

Assume that not all pairs $pt$ with $p,t\in B$ have the same value for $d_0(pt)$. We can choose $ij$ with $i,j\in B$ so that at least $\frac{1}{2}\binom{|B|}{2}$ pairs $pt$ of $\binom{B}{2}$ have $d_0(ij)\ne d_0(pt)$. Then the number of pairs $pt$ of $\binom{B}{2}$ in $\binom{B\backslash\{i,j\}}{2}$ with $d_0(ij)\ne d_0(pt)$ is at least 
$\frac{1}{2}\binom{|B|}{2}-2|B|$.  Now using the second inequality of (\ref{clique1100}) with $|{\cal U}(ij)|+|{\cal T}(ij)|\le (c_3+c_4)m$ and $|B|\le m$,  we can find a pair
$k,l\in B\backslash\{i,j\}$ with $d_0(ij)\ne d_0(kl)$, $kl\notin {\cal U}(ij)\cup{\cal T}(ij)$. 
By definition of $B$, 
$$d_0(ij)+d_1(ij)=q-3,\quad d_0(kl)+d_1(kl)=q-3.$$
We may assume without loss of generality that
$d_0(kl)<d_0(ij)$, $d_1(kl)>d_1(ij)$ and then
$$d_0(ij)+d_1(kl)\ge q-2$$
We also have $\epsilon(ij)=\epsilon(kl)=0$. Then $A$ has a column of column sum 2 and a  column of column sum $m-2$ both with 1's in rows $k,l$ and 0's in rows $i,j$. Also we have $d_0(ij)+d_1(kl)$ columns with 1's in rows $k,l$ and 0's in rows $i,j$ since $kl\notin {\cal U}(ij)\cup{\cal T}(ij)$. But then $A$ has $q\cdot (\1_2\0_2)$, a contradiction.
We conclude that all pairs $pt$ with $p,t\in B$ have the same value for $d_0(pt)$.

We follow our proof of Lemma~\ref{onlybd1100} using (\ref{a2a3a4eq1100}) and deduce that $E_0\cup E_1=\emptyset$ and so $a_2=2\binom{m}{2}$. Also we deduce that 
$$\sum_{ij}t(ij)=3a_4$$
and as a result we can deduce that any column $\alpha$ in $A_4$ either has column sum 4 with exactly one 1 in a row of $B$ or has column sum $m-4$ with exactly one 0 in a row of $B$. 

Assume $a_4>0$ and consider $\alpha$ in $A_4$, say with column sum 4 and with 1's in rows $i,j,k,l$ where $i\in B$ and 
$j,k,l\in\{1,2,\ldots ,m\}\backslash B$. Choose $r,s\in B\backslash i$ so that $d_0(rs)+d_1(rs)=q-3$ and with $rs\notin {\cal T}(ij)\cup {\cal U}(ij)$ (using first inequality of (\ref{clique1100})).  Column $\alpha$ has 1's in rows $i,j$ and 0's in row $r,s$. Using $E_0\cup E_1=\emptyset$, we deduce that
$d_1(ij)+d_0(rs)\le q-3-1$ and
$d_0(ij)+d_1(rs)\le q-3$  else if either inequality  is violated we create $q\cdot(\1_2\0_2)$. We deduce
$d_0(ij)+d_1(ij)\le (q-3)-1$. This yields a slight variant of (\ref{sumd0d1}): 
$$\sum_{ij}(d_0(ij)+d_1(ij)) \leq 
\sum_{ij}{(q-3+\epsilon(ij)-t(ij))}-1.$$
The extra `-1' is sufficient to obtain a contradiction when we substitute for $a_2,a_3,a_4$ in (\ref{a2a3a4eq1100}). We then deduce $a_4=0$.

With $a_4=0$ and $a_2=2\binom{m}{2}$, we deduce
$a_3=\frac{q-3}{3}\binom{m}{2}$ using (\ref{a2a3a4eq1100}). Given that $\epsilon(ij)=0$ for all $ij$ and using Lemma~\ref{dupperbd1100}, we deduce 
$d_0(ij)+d_1(ij)=q-3$ for all pairs $ij$ and so $B=\{1,2,\ldots ,m\}$. Our above arguments tell us $d_0(pt)$ is the same for every choice $p,t\in B$, allowing us to conclude that there exists positive integers $a,b$, $a+b=q-3$ so that for all pairs
$ij$, $d_0(ij)=a$ and $d_1(ij)=b$. From this we can readily conclude that the columns of column sum 3 correspond to a 2-design $S_{a}(2,3,m)$ and the columns of column sum $m-3$ correspond to the (0,1)-complement of a 2-design $S_{b}(2,3,m)$. \qed

\end{document}